\documentclass[11pt]{article}
\usepackage{graphicx}
\usepackage{amssymb}
\date{}
\newtheorem{teo}{Theorem}[section]
\newtheorem{cor}[teo]{Corollary}
\newtheorem{defi}[teo]{Definition}
\newtheorem{lema}[teo]{Lemma}
\newtheorem{prop}[teo]{Proposition}
\newtheorem{rem}[teo]{Remark}

\newcommand{\pf}{{\bf Proof.\ \ }}

\newcommand{\text}{}

\newcommand{\Ad}{\mathop{\rm Ad}}

\newcommand{\Tr}{\mathop{\rm Tr}}

\newcommand{\sig}{\mathop{\rm sig}}

\title{A METHOD FOR THE
RESOLUTION OF THE JACOBI EQUATION $Y'' + RY = 0$ ON THE MANIFOLD
$Sp(2)/SU(2)$.}
\author{A. M. Naveira and A. Tarr\'{\i}o {\footnote{Work
partially supported by a Research Project MTM 2004-06015-C02-01 
(first author)
and by a Research Project PGIDIT05PXIB16601PR (second author).
}}}

\begin{document}

\maketitle

\begin{abstract}	In this paper a method
for the resolution of the differential equation of the Jacobi vector 
fields in the
manifold $V_1 = Sp(2)/SU(2)$ is exposed. These results are applied 
to determine
areas and volumes of geodesic spheres and balls.
\end{abstract}



\section*{Introduction}	The resolution
of the Jacobi equation on a Riemannian manifold can be quite a 
difficult
task. In the Euclidean space the solution is trivial. For the 
symmetric
spaces, the problem is reduced to a system of differential equations 
with
constant coefficients. In the specialized bibliography, particularly 
in
\cite{[Gr]}, the explicit solutions of these systems are found as 
well as their
application to the determination of areas and volumes. In  \cite 
{[Ch1],[Ch2]} a
partial solution of this problem for the manifolds $V_1 = 
Sp(2)/SU(2)$ and $V_2 
=
SU(5)/{Sp(2)\times S^1}$ is obtained by I. Chavel. It is well known that these 
manifolds are
nonsymmetric normal homogeneous spaces of rank 1 \cite [p.237]{[B]}. 
The manifold
$V_1$ appears in \cite{[B]} and in the book of A. L. Besse \cite [p.203]{[Bs]} as an
exceptional naturally reductive homogeneous space. For  naturally
reductive compact homogeneous spaces, Ziller \cite{[Z]} solves the 
Jacobi equation  
working with the canonical connection, which is natural for the 
nonsymmetric
naturally reductive homogeneous spaces; but the solution can be 
considered
of qualitative type, in the sense that it does not allow to obtain 
in an easy
way the explicit solutions of the Jacobi fields for any particular 
example and
for an arbitrary direction of the geodesic. The method used by 
Chavel, which
allows him to solve the Jacobi equation in some particular 
directions, is
based on the use of the canonical connection. Nevertheless, his method
 does not seem to  apply in
a simple way to the resolution of  the Jacobi
equation along a unit geodesic of arbitrary direction. In 
\cite 
{[Ch1],[Ch2]} the same author shows the existence of anisotropic Jacobi fields;  
that is, they do 
not 
come from geodesic variations in the isotropic subgroup. Also, the 
Jacobi equations on
a Riemannian manifold appear in a natural way in the theory of 
Fanning curves 
\cite{[AP-D]}.

In this
paper, always working with the Levi-Civita connection and using an interesting geometric
result of Tsukada \cite{[T]}, the Jacobi equation along a 
unit 
geodesic of
arbitrary direction is solved. Also, the solutions are applied to 
obtain the
area of the geodesic sphere and the volume of the geodesic ball of 
radius $t$ in
the manifold $V_1= Sp(2)/SU(2)$.	In \S 1, for an arbitrary Riemannian manifold, using the
induction method,  a 
recurrent formula
for the $i$-th covariant derivative of the Jacobi operator
$R_t = R(\cdot ,\gamma')\gamma' (t)$ along the geodesic $\gamma$ is given. 	In \S 
2, using
the result of the previous section, the expression of the covariant 
derivative
of the curvature tensor at the point $\gamma(0)$ is obtained for an arbitrary naturally 
reductive
homogeneous spaces $M =G/H$, in terms of the
brackets of the Lie algebra of $G$. In order to  obtain this result the 
induction
method is used again. 	
In the following sections, the previous 
results are
applied to the normal homogeneous space $V_1$. So, in 
\S 3, always
working with a unit geodesic $\gamma$ of arbitrary direction, 
the values at $\gamma(0)$ of the Jacobi tensor $R_0$, and its covariant derivatives
$R\sp{1)}_0$ and
$R\sp{2)}_0$  are
determined. In Lemma 3.1 it is 
proved that, for a unit geodesic
$\gamma$,
$$R\sp{3)}_0= -{\vert \vert \gamma'\vert \vert }^2 R\sp{1)}_0 = -
R\sp{1)}_0$$		
and
 $$ R\sp{4)}_0= -{\vert \vert \gamma'\vert \vert}^2 
R\sp{2)}_0
= - R\sp{2)}_0.$$ 	
This section ends by proving
that
$$R\sp{2n)}_0 = (-1)\sp{n-1} R\sp{2)}_0$$
	and
$$R\sp{2n+1)}_0 = (-1)\sp{n} R\sp{1)}_0.$$
Using
the Taylor development, at the point $\gamma(0)$, of the Jacobi 
operator, it is
possible to obtain quite a simple expression of the Jacobi operator 
$R_t$
along the geodesic, as well as that of its derivatives.  In fact, the
explicit expression for $R_t$ is 
$$R_t= R_0 + R\sp{2)}_0 + R\sp{1)}_0\sin t - 
R\sp{2)}_0
\cos t.$$			
It seems interesting to remark that while
 D'Atri and
Nickerson \cite{[DN1], [DN2]} impose conditions on the derivatives of the
Jacobi operator, known as Ledger's conditions of odd order, in our case,
conditions are imposed by the geometric properties of the manifold.

 In \S 4 the
Jacobi equation with predetermined initial values is solved and the 
formal
expressions of the area and the volume of the geodesic sphere and 
the ball of
radius $t$ are obtained.	In a forthcoming paper the problem of 
determining
the areas of tubular hypersurfaces and the volumes of tubes around 
compatible
submanifolds will be approached. Given its generality, we hope 
that this
method could also be applied to solve the Jacobi equation in 
several other examples of naturally reductive homogeneous spaces.

\section{A formula for the
covariant derivative of the Jacobi operator in a Riemannian 
manifold.}	Let $M$
be an $n$-dimensional, connected, real analytical Riemannian 
manifold, $g = < , >$ its Riemannian metric, $m\in M$, $v\in T_m M$ a unit tangent 
vector and 	$\gamma : J \to M$ a geodesic in $ M$ defined on some 
open
interval $J$ of ${\mathbb {R}}$ with $0   \in J$, $m=\gamma(0)$. 
For a geodesic $\gamma(t)$ in $M$ the
associated Jacobi  operator
$R_t$ is the self-adjoint tensor field along $\gamma$ defined by 
$$R_t : =
R(\cdot ,\gamma')\gamma' (t)$$ for the curvature tensor we follow the 
notations of \cite{[KN]} . The covariant derivative
$R\sp{i)}_t$ of the Jacobi  operator
$R_t$ along $\gamma$ is the self-adjoint tensor field 
defined
by 
$$R\sp{i)}_{t}:= (\nabla\sb{\gamma'} \stackrel{i)} \cdots 
\nabla\sb{\gamma'}R)(\cdot,\gamma')\gamma' (t),$$	
where $\nabla$ is the Levi-Civita
connection associated to the metric. Its value at $\gamma (0)$ will be denoted by
$$R\sp{i)}_{0}:= (\nabla\sb{\gamma'} \stackrel{i)}\cdots  
\nabla\sb{\gamma'}R)(\cdot,\gamma')\gamma'(0)$$

and we denote $R\sp{_{0)}}_t= R_t.$	

First, we
prove two combinatorial lemmas for later use.

\begin{lema}  For $i
\leq 2k$ we have:
\begin{enumerate}
\item[a)] $$\left( {\begin{array}{*{20}c}   {2k + 2}  \\   i  \\
\end{array} } \right) = \left( {\begin{array}{*{20}c}   {2k}  \\   i  
\\
\end{array} } \right) + 2\left( {\begin{array}{*{20}c}   {2k}  \\   
{i - 1} 
\\ \end{array} } \right) + \left( {\begin{array}{*{20}c}   {2k}  \\   
{i -
2}  \\\end{array} } \right);$$						 		
\item[b)]	$$\left(
{\begin{array}{*{20}c}   {{{2k + 2}}}  \\   {{{2k + 1}}}  \\ 
\end{array} }
\right) = \left( {\begin{array}{*{20}c}   {{{2k}}}  \\   {{{2k - 
1}}}  \\
\end{array} } \right) + 2;$$		
\item[c)] $$\left( {\begin{array}{*{20}c}   
{{{2k +
2}}}  \\   {{{2k + 2}}}  \\ \end{array} } \right) = \left(
{\begin{array}{*{20}c}   {{{2k}}}  \\   {{{2k}}}  \\ \end{array} } 
\right) =
1.$$
\end{enumerate}
\end{lema} 	
The proof is a trivial consequence of some properties of the
combinatorial numbers.
\begin{lema} 
\[\sum\limits_{j = 1}^i {( - 1)^j \left(
{\begin{array}{*{20}c}   {k + 1}  \\   j  \\ \end{array} } \right)} 
\left(
{\begin{array}{*{20}c}   {k - j + 1}  \\   {i - j}  \\ \end{array} }
\right) = \left(
{\begin{array}{*{20}c}   {k + 1}  \\   i  \\ \end{array} } \right).\]
\end{lema}
The proof follows at once by using the formula
\[\left(
{\begin{array}{*{20}c}   { - x}  \\   n  \\ \end{array} } \right) = 
( - 1)^n
\left( {\begin{array}{*{20}c}   {x + n - 1}  \\   n  \\ \end{array} 
}
\right)\]	where $x\in {\mathbb {Z}}$, and also the Vandermonde's 
identity
\[\left(
{\begin{array}{*{20}c}   {x + y}  \\   n  \\ \end{array} } \right) =
\sum\limits_{j = 0}^n {\left( {\begin{array}{*{20}c}   x  \\   j  \\
\end{array} } \right)} \left( {\begin{array}{*{20}c}   y  \\   {n - 
j}  \\
\end{array} } \right)\]
with $x, y \in {\mathbb {Z}}.$
\begin{teo}
For $n \geq 1$ we have $$\nabla_{\gamma'}	\stackrel{
n)}\cdots\nabla_{\gamma'}R(X, \gamma') \gamma' =  \sum\limits_{i = 0}^n 
{\left(
\begin{array}{*{20}c}   n  \\   i  \\ \end{array}  \right)} 
R\sp{\sp{n-i)}}\sb{t}(\nabla_{\gamma'} \stackrel{
i)}\cdots\nabla_{\gamma'}X).$$
\end{teo}
\pf We prove this by induction. For $n
= 1$, we have
\begin {equation}\label {(1.3)}
\nabla_{\gamma'}R(X, \gamma') \gamma' =
(\nabla_{\gamma'}R)(X,\gamma')\gamma'+ R(\nabla_{\gamma'}X, 
\gamma')\gamma'    
\end {equation} 
that is $$\nabla_{\gamma'}R(X, \gamma') \gamma' = R 
\sp{1)}\sb{t}(X)+ R
\sp{_0)}\sb{t}(\nabla_{\gamma'}X)$$ and so the result is true for $n = 
1$. 
Next,
suppose that Theorem 1.3 holds for $n= k$. 
Then we have
$$\nabla_{\gamma'} 	\stackrel{
k)}\cdots\nabla_{\gamma'}R(X, \gamma') \gamma' =  \sum\limits_{i = 
0}^k
{\left( {\begin{array}{*{20}c}   k  \\   i  \\ \end{array} } 
\right)}
R\sp{k-i)}_t(\nabla_{\gamma'}	\stackrel{ i)}\cdots  
\nabla_{\gamma'}X). $$	
Taking the covariant derivative, we obtain 
\begin{eqnarray*}\nabla_{\gamma'}(\nabla_{\gamma'}
\stackrel{ k)}\cdots 
\nabla_{\gamma'}R(X,\gamma')\gamma')&=&\nabla_{\gamma'}
\stackrel{k+1)}\cdots\nabla_{\gamma'}R(X, \gamma') \gamma'
\\&=&\nabla_{\gamma'}( \sum\limits_{i = 0}^k \left( 
\begin{array}{*{20}c}   k
 \\   i  \\\end{array}  \right)
R\sp{k-i)}_t(\nabla_{\gamma'}	\stackrel{ i)}\dots  
\nabla_{\gamma'}X)).\\\end{eqnarray*}
By applying (\ref {(1.3)}) to each term, it is possible
to write
\begin{eqnarray*}&&\nabla_{\gamma'}
\stackrel{k+1)}\cdots\nabla_{\gamma'}R(X, \gamma') \gamma' \\&=&
\sum\limits_{i = 0}^{k } {\left( {\begin{array}{*{20}c}   k  \\   i  
\\
\end{array} } \right)}  [R\sp{k+1-i)}_t(\nabla_{\gamma'}
\stackrel{i)}\cdots\nabla_{\gamma'}X)  +
R\sp{k-i)}_t(\nabla_{\gamma'}
\stackrel{i+1)}\cdots\nabla_{\gamma'}X)]  \\&=&\left( 
\begin{array}{*{20}c}   {k
+ 1}  \\   0  \\\end{array}  \right)R\sp{k+1)}_t(X) + 
\sum\limits_{i =
0}^{k - 1}[\left( \begin{array}{*{20}c}k \\i  \end{array} 
\right) + \left( \begin{array}{*{20}c}k  \\i+1  \\\end{array} 
\right)]
R\sp{k-i)}_t(\nabla_{\gamma'}
\stackrel{i+1)}\dots\nabla_{\gamma'}X) \\&&+ \left( 
\begin{array}{*{20}c}k 
\\k  \\\end{array}  \right)R\sp{0)}\sb{t}(\nabla_{\gamma'} 
\stackrel{ k+1)}\dots\nabla_{\gamma'}X).\\\end{eqnarray*}	
Now, by applying basic properties of combinatorial numbers
we have
$$\nabla_{\gamma'} 	\stackrel{
k+1)}\cdots\nabla_{\gamma'}R(X, \gamma') \gamma' = \sum\limits_{i = 
0}^{k+1
}{\left( {\begin{array}{*{20}c}   {k+1}  \\   i  \\ \end{array} } 
\right)}
R\sp{k+1-i)}_t(\nabla_{\gamma'}	\stackrel{ i)}\cdots  
\nabla_{\gamma'}X) $$	
and the result follows.

\begin{cor} We have $$R\sp{\sp{n)}}\sb{t}(X)=\nabla_{\gamma'}
\stackrel{n)}\cdots\nabla_{\gamma'}R(X, \gamma') \gamma' - 
\sum\limits_{i =
1}^{n} {\left(
\begin{array}{*{20}c}   n  \\   i  \\ \end{array}  \right)} R\sp{n-
i)}_t
(\nabla_{\gamma'} \stackrel{ i)}\cdots 
\nabla_{\gamma'}X).$$
\end{cor}

\section{An algebraic expression for the
covariant derivative of the Jacobi operator on a naturally 
reductive  
homogeneous space} Let	$G$ be a Lie group, $H$ a closed subgroup, 
$G/H$ the
space of left cosets of $H$, $\pi : G\to G/H$ the natural 
projection.	For $r \in G$ we denote by $\tau$ the induced action of $ G$ 
on $G/H $
given by $\tau(r)(sH) = rsH$, $r, s \in G$. The Lie algebras of $G$ 
and $H$
will be denoted by ${{\bf g}}$ and ${\bf h}$, respectively and ${\bf 
m} ={{\bf
g}}/{\bf h}$ is a vector space which we identify with the tangent 
space to $G/H$
at  $o = \pi(H)$. 	An affine connection on $G/H$ is said to be invariant if 
it is
invariant under $\tau(r)$ for all $r \in G$. 	It is well known that 
it is
possible to define in a natural way on ${{\bf g}}$ an $\Ad$-invariant 
metric by 
$<u, v> = \Tr(uv^{t})$, $u, v \in {{\bf g}}$. Let 
$\nabla$ be the associated Levi-Civita connection. It is
well-known \cite [Ch.X, p.186]{[KN]} that there exists an invariant affine 
connection $D$ 
on
$G/H$ (the {\it  canonical connection}) whose torsion $T$ and 
curvature $B $ 
tensors
are also invariant. In the following we always work with $\nabla $.
\begin{defi} {\em \cite [p.202] {[Ch2], [KN]}} $M = G/H$ is said 
to be a
\begin{enumerate}
\item[(a)]
{ {\em
Reductive homogeneous space}} if the Lie algebra ${{\bf g}}$ admits 
a vector
space decomposition ${{\bf g}} = {\bf h} + {\bf m}$ such that $[{\bf 
m}, {\bf
h}] \subset {\bf m}$. In this case ${{\bf m}}$ is identified with 
the tangent
space at the origin $o = \pi(H)$.
\item[(b)] {{\em Riemannian homogeneous space}} if $G/H $ is
a Riemannian manifold such that the metric is preserved by $\tau(r)$ 
for all  
$r \in G$.

\item[(c)] {{\em Naturally reductive  Riemannian homogeneous 
space}} if 
$G/H$,
with a H-invariant Riemannian metric, admits an $\Ad(H)$-invariant 
decomposition 
 ${{\bf g}}= {\bf h}+ {\bf m}$ satisfying the condition
$$<[u, v]_{\bf m} , w> + <v, [u, w]_{\bf m}> = 0$$ for $ u, v, w \in 
{\bf m}$.

\item[(d)] {{\em Normal Riemannian homogeneous space}} if the 
metric
on $G/H$ is obtained as follows: there exists a positive definite 
inner
product $< , >$ on ${{\bf g}}$ satisfying $$<[u, v], w> = <u, [v, 
w]>$$ for all
$u, v, w  \in {{\bf g}}$. Let ${\bf m} = {{\bf g}}/{\bf h}$ be 
the
orthogonal complement of {\bf h}. Then the decomposition $({{\bf 
g}}, {\bf h})$
is reductive, and the restriction of the inner product to ${\bf m}$ 
 induces a Riemannian 
metric
on $G/H$, referred to as {{\em normal}}, by the action of $G$ on
$G/H$.
\end{enumerate}
\end{defi}	
From now on we will assume that $G/H$ is a naturally reductive space. If we define $\Lambda
\colon  {\bf m}
\times {\bf m}\to {\bf m}$ by $$\Lambda(u)v = (1/2)[u, v]_{\bf m}$$ for
$u, v\in {\bf m}$, we can identify $\nabla$ and $\Lambda$. Evidently,
$\Lambda(u)$ is a skew-symmetric linear endomorphism of $({\bf m}, < 
, >)$.
Therefore $e^{\Lambda(u)}$ is a linear isometry of $({\bf m}, < , 
>)$. Since
the Riemannian connection is a natural torsion free connection on 
$G/H$, we have \cite [Vol. II, Ch.X]{[Tj],[KN]}:
\begin{prop}
The following properties hold:
\begin{enumerate}
\item[(i)] 
For each $v
\in{\bf{ m}}$, the curve $\gamma(t) = \tau(\exp tv)(o)$ is a 
geodesic with
$\gamma(0) = o$, $\gamma'(0) = v$.
\item[(ii)] 
The parallel translation along
$\gamma$ is given as follows: $$\tau(\exp tv)_{*} e\sp{-
t\Lambda(v)}\colon T_oM
\to T_{\gamma(t)}{M}.$$
\item[(iii)]
 The (1,3)-tensor $R_t$ on ${\bf m}$ obtained
by the parallel translation of the Jacobi operator along $\gamma$ is 
given as
follows:
$$R_t = e \sp{t\Lambda(v)} R_0.$$
\end{enumerate}
\end{prop}
Above, $R_0$
denotes the Jacobi operator at the origin $o$ and 
$e\sp{t\Lambda(v)}$
denotes the action of $e\sp{t\Lambda(v)}$ on the space $R({\bf m})$ 
of
curvature tensors on ${\bf m}$.

\begin{prop} {\em \cite [Vol. II, p. 202]{[Ch2],[KN]}} Let $\gamma(t)$ be a 
geodesic with 
$\gamma(0) = o$, for
$v = \gamma'(0) \in {\bf m}$. If X is a  differentiable vector field 
along
$\gamma$, then$$R_0(X)=  -[[X, v]_{\bf h},v] - (1/4)[[X,v]_{\bf 
m},v]_{\bf
m}.$$
\end{prop}
\begin{prop}Under the same hypothesis that in Proposition
2.3, we have, for $ n>0$,
\begin {equation}\label{(2.1)}
(-1)\sp{n-1}  
2\sp{n}R\sp{n)}_0(X)=\sum\limits_{i
= 0}^{n} {( - 1)^i \left( \begin{array}{*{20}c}n  \\i  \\\end{array} 
\right) [[[X,v]_{\bf m}, \dots , v   ]^{i+1)}_{\bf 
h}  , \dots , v]_{\bf m}    } \end {equation}

 where for each term of the sum we have $n+2$ brackets and the exponent $i+1)$ means the 
position of
the bracket valued on ${\bf h}$.\end{prop}
\pf Using Proposition 2.3,
Corollary 1.4 and the fact that 
\begin{equation} \label{(2.2)}
\nabla_{ X}Y = (1/2)[X,Y]_{\bf m},
\quad  X,Y \in 
{\bf m}, 
\end{equation}			
we have immediately that (\ref {(2.1)}) is 
verified for 
$n = 1$.
Next, suppose that this formula holds for $n= k$; then
\begin{equation} \label{(2.3)}
(-1)^{k-1}2^{k}
R\sp{k)}_0 (X)=\sum\limits_{{{i = 0}}}^{{k}} {( - 1)^i \left(
{\begin{array}{*{20}c}   {{k}}  \\   {{i}}  \\\end{array} } \right)} 
[[[X,v]_{\bf m}, \dots , v]^{i+1)}_{\bf h}, \dots ,
v]_{\bf m}.   \end {equation} 	
Using now Corollary 1.4, we
have$$R\sp{k+1)}_0 (X)=  \nabla_v \stackrel{
k+1)}\cdots\nabla_vR(X,v)v -\sum\limits_{{{i = 1}}}^{{k+1}} { \left(
{\begin{array}{*{20}c}   {{k+1}}  \\   {{i}}  \\\end{array} } 
\right)} 
R\sp{k+1-i)}_0(\nabla_v \stackrel{ i)}\cdots \nabla_vX).$$  In each
term we take into account Proposition 2.3, formulae (\ref{(2.2)}) and (\ref{(2.3)}), so we
obtain 
\begin{eqnarray*}R\sp{k+1)}_0(X)&=& (-1)\sp{k}  
\frac{1}{{2^{k
+ 1} }}[[X,v]_{\bf h}, \dots ,v]_{\bf m}
\\&&-\sum\limits_{i = 1}^{k + 1} {\left( {\begin{array}{*{20}c}   {k 
+ 1} 
\\   i  \\\end{array} } \right)}  (-1)\sp{k-1}\frac{1}{{2^{k + 1 - 
i}
}}\\&& \left((-1)^i \frac{1}{{2^i }}\sum\limits_{j = 0}^{k + 1 - i} 
{\left(
{\begin{array}{*{20}c}   {k + 1 - i}  \\   j  \\\end{array} } 
\right)} 
(-1)^j [[[X, v]_{\bf m},\dots , v]_ {\bf h}^{i+j+1)}, 
\dots , v]_{\bf m}\right)\end{eqnarray*} Let us remark that the sum 
of the terms with all
brackets  estimated
in {\bf m} is $0$. By the other hand the terms that have the bracket estimated in 
${\bf h}$ in
the $(i+1)$-position are
\begin{eqnarray*}
&&-\left( 
{\begin{array}{*{20}c}  
{k + 1}  \\   1  \\ \end{array} } \right)\frac{1}{{2^k }}(-1)^{k-1}
(-1)^{i-1}\frac{1}{2}( - 1) \left( {\begin{array}{*{20}c}   k  \\ 
{i - 1}
 \\ \end{array} } \right)\\
 &&\\%
 &&-\left( {\begin{array}{*{20}c}   {k + 1}  \\  
2  \\ \end{array} } \right)\frac{1}{{2^{k - 1} }}\frac{1}{{2^2 }} (-
1)^{k-2}
(-1)^{ i-2}\left( {\begin{array}{*{20}c}   {k - 1}  \\   {i - 2} 
\\ \end{array} } \right)\\
 &&\\%
&& - \cdots  -\left( {\begin{array}{*{20}c}   
{k +
1}  \\   i  \\\end{array} } \right)\frac{1}{{2^{k + 1 - i} }} (-
1)^{k-i}
(-1)^i \frac{1}{{2^i }}\left( {\begin{array}{*{20}c}   {k + 1 - i}  
\\   0 
\\\end{array} } \right) \\
&&\\&&\\%
&= &( - 1)^k \frac{1}{{2^{k + 1} }}( - 
1)^i
(\left( {\begin{array}{*{20}c}   {k + 1}  \\   1  \\\end{array} }
\right)\left( {\begin{array}{*{20}c}   k  \\   {i - 1}  
\\\end{array} }
\right)\\
&&\\%
&&  - \left( {\begin{array}{*{20}c}   {k + 1}  \\   2 
\\\end{array} } \right)\left( {\begin{array}{*{20}c}   {k - 1}  \\   
{i - 2}
 \\\end{array} } \right)\\
&&\\ %
 && +\cdots + ( - 1)^{i - 1} \left(
{\begin{array}{*{20}c}   {k + 1}  \\   i  \\\end{array} } 
\right)\left(
{\begin{array}{*{20}c}   {k + 1 - i}  \\   0  \\\end{array} }
\right)).
\end{eqnarray*} 

Using Lemma 1.2, the last expression equals	
$$
( - 1)^k
\frac{1}{{2^{k + 1} }}( - 1)^i  \left( {\begin{array}{*{20}c}   {k + 
1}  \\ 
 i  \\\end{array} } \right)
$$
 
 so the formula (\ref{(2.1)}) is true for $n=k+1$ and 
this finishes the proof.		

\section{ An
explicit form for the Jacobi operator on the 
manifold $V_1 = Sp(2)/SU(2)$.}	
We consider the Lie group $Sp(2)$ and the
subgroup $SU(2)$. It is well known that $V_1 = Sp(2)/SU(2)$ is a 
normal
naturally reductive Riemannian homogeneous space \cite{[B],[Ch1],[Ch2]}. 
We denote
by $sp(2)$ and $su(2)$ the Lie algebras of $Sp(2)$ 
and $SU(2)$
respectively. Using the notations of \cite{[Ch2]} it is known that 
an element of
the Lie algebra $sp(2)$ is a skew-Hermitian matrix of the form
\[\left(
{\begin{array}{*{20}c}   {a_{11} } & {a_{12} } & {a_{13} } & {a_{14} 
}  \\  
{ - \overline a _{12} } & { - a_{11} } & {\overline a _{14} } & { - 
\overline a
_{13} }  \\   { - \overline a _{13} } & { - a_{14} } & {a_{33} } & 
{a_{34} } 
\\   { - \overline a _{14} } & {a_{13} } & { - \overline a _{34} } & 
{ - a _{33} }  \\ \end{array} } \right)\]
where
$a\sb{11}$, $a\sb{33}$ are pure imaginary numbers and the other $a _{ij}$ are arbitrary 
complex
numbers. Let $S_{i}$,  $i = 1,\dots, 10$ be the matrices of ${sp(2)}$ 
such
that\begin{eqnarray*}	S_1&: &a\sb{11}= - a\sb{22} = i;\\
S_2&: &a\sb{33}= - a\sb{44} = i;
\\		S_3 &: &a\sb{12}= - a\sb{21} = 1; 
\\			S_4 &: &a\sb{12}= 
a\sb{21} = i;
\\			S_5 &:& a\sb{34}= - a\sb{43} = 1;
 \\		S_6 &:& a\sb{34}= 
a\sb{43} = i	;
\\		S_7 &: &a\sb{13} = - a\sb{31} = a\sb{24} = - a\sb{42} =
1;
 \\	S_8 &: &a\sb{13} =  a\sb{31} =  a\sb{24} =  a\sb{42} = i;
\\	S_{9}& :&
a\sb{14} = - a\sb{41} = a\sb{23} = - a\sb{32} = 1;
\\	S_{10} &: & a\sb{14} = 
a\sb{41}=  a\sb{23} = a\sb{32} = i;
\\	\end{eqnarray*} 	
	the other $a_{ij}$ being zero in all cases. Evidently $\{S_i\}$ is 
an adapted basis
of $sp(2)$. We construct another basis $\{Q_i\}$ as follows:	
\small
$$\left(
{\begin{array}{*{20}c}   {Q_1 }  \\   {Q_2 }  \\   {Q_3 }  \\   {Q_4 
}  \\
  {Q_5 }  \\   {Q_6 }  \\   {Q_7 }  \\   {Q_8 }  \\   {Q_9 }  \\   
{Q_{10}
}  \\\end{array} } \right)	= 	\left( {\begin{array}{*{20}c}   
{1/2} & { -
3/2} & 0 & 0 & 0 & 0 & 0 & 0 & 0 & 0 
 \\   0 & 0 & {\sqrt {5/2} } & 0 
& 0 & 0 &
0 & 0 & 0 & 0  
\\   0 & 0 & 0 & {\sqrt {5/2} } & 0 & 0 & 0 & 0 & 0 & 
0 
 \\  
0 & 0 & 0 & 0 & {\sqrt 6 /2} & 0 & { - \sqrt 2 /2} & 0 & 0 & 0  
\\   
0 & 0 & 0
& 0 & 0 & {\sqrt 6 /2} & 0 &{ - \sqrt 2 /2} & 0 & 0  \\   0 & 0 & 0 & 0 & 0 & 0 & 
0 & 0 &
{\sqrt {5}/2 } & 0  \\   0 & 0 & 0 & 0 & 0 & 0 & 0 & 0 & 0 & {\sqrt 
{5}/2 } 
\\   {3/2} & {1/2} & 0 & 0 & 0 & 0 & 0 &0& 0 & 0  \\   
0 & 0
& 0 & 0 & 1 & 0 & {\sqrt {3}/2 } & 0 & 0 & 0  \\   0 & 
0 & 0 & 0
& 0 & 1 & 0 & {\sqrt {3}/2 } & 0 & 0 \\ \end{array} }
\right)\left( {\begin{array}{*{20}c}   {S_1 }  \\   {S_2 }  \\   
{S_3 }  \\
  {S_4 }  \\   {S_5 }  \\   {S_6 }  \\   {S_7 }  \\   {S_8 }  \\   
{S_9 } 
\\   {S_{10} }  \\ \end{array} } \right).$$
\normalsize	
We have \cite{[Ch2]}:
\begin{enumerate}
\item[i)] 	

 If for an inner product on $sp(2)$ we take $<A, B> = -(1/5) 
\Tr(AB)$,
then $\{Q_{1}, \dots, Q_{10}\}$ is an orthonormal basis of $sp(2)$.

\item[ii)]The inner product is invariant under $\Ad(Sp(2))$.		
					
\item[iii)] Finally, one can show that ${\bf h} =$ linear span of $ 
\{Q_{8},Q_{9},Q_{10}\} $
is Lie diffeomorphic to $su(2)$ and therefore the group generated by 
{\bf h} is
analytically isomorphic to $SU(2)$.
\end{enumerate}
								 
The previous decomposition is taken
from \cite [p.234]{[B]}. If we call $ {\bf m} = sp(2)/su(2)$, we 
know that $\{Q_{1},\dots, Q_{7}\}$ is an
adapted basis for ${\bf m}$. It is 
immediate to
prove that the brackets are given by the following relations 
\begin{equation}
 \begin{array} {lcl} 
 \left[Q_{1}, Q_{2}\right] = Q_{3},\hfill	&&		\left[Q_{1}, Q_{3}\right] = - Q_{2},\\ 
\left[Q_{1}, Q_{4}\right] = - Q_{5} - \sqrt{6}Q_{10}, && \left[Q_{1},  Q_{5}\right] =  Q_{4} + \sqrt{6}Q_{9},\\ 
\left[Q_{1}, Q_{6}\right] = -Q_{7},		&&  \left[Q_{1}, Q_{7}\right] = Q_{6},\\ 
\left[Q_{1}, Q_{8}\right] = 0,		&&	\left[Q_{1}, Q_{9}\right] = - \sqrt{6}Q_{5},\\ 
\left[Q_{1}, Q_{10}\right] = \sqrt{6}Q_{4},	&&	\left[Q_{2}, Q_{3}\right] = Q_{1} + 3  Q_{8},\\ 
\left[Q_{2},Q_{4}\right] = Q_{6}, 	&&		\left[Q_{2}, Q_{5}\right] = - Q_{7},\\ 
\left[Q_{2}, Q_{6}\right] = - Q_{4}+\sqrt{3/2}Q_{9},&& \left[Q_{2}, Q_{7}\right] = Q_{5} - \sqrt{ 3/2} Q_{10},\\ 
\left[Q_{2}, Q_{8}\right] = - 3Q_{3}, 	&&	\left[Q_{2}, Q_{9}\right] = - \sqrt{ 3/2} Q_{6},\\ \label{CORCHETE} 
\left[Q_{2}, Q_{10}\right] = \sqrt{ 3/2}Q_{7},	&&	\left[Q_{3}, Q_{4}\right] = Q_{7}, \\ 
\left[Q_{3}, Q_{5}\right] =Q_{6}, 		&&	\left[Q_{3},Q_{6}\right] = -(\sqrt{  2}  /2)( \sqrt{  2}  Q_{5} -\sqrt{  3}   Q_{10}),\\ 
\left[Q_{3},Q_{7}\right]=-(\sqrt{  2}  /2)( \sqrt{  2} Q_{4}-\sqrt{  3}  Q_{9}),&& \left[Q_{3}, Q_{8}\right] = 3 Q_{2},\\ 
\left[Q_{3}, Q_{9}\right] = -\sqrt{ 3/2}  Q_{7}, &&	\left[Q_{3}, Q_{10}\right] = -\sqrt{ 3/2}  Q_{6},\\ 
\left[Q_{4}, Q_{5}\right] = - Q_{1} + Q_{8},&&	\left[Q_{4}, Q_{6}\right] = Q_{2} + \sqrt{ 5/2}  Q_{9},\\ 
\left[Q_{4}, Q_{7}\right] =Q_{3} +\sqrt{ 5/2}   Q_{10},&&	\left[Q_{4}, Q_{8}\right] = - Q_{5},\\ 
\left[Q_{4}, Q_{9}\right] = - \sqrt{ 5/2}   Q_{6},	&&\left[Q_{4}, Q_{10}\right] = - 2\sqrt{ 3/2}   Q_{1} - \sqrt{ 5/2}   Q_{7},\\ 
\left[Q_{5}, Q_{6}\right] =Q_{3} -\sqrt{ 5/2}  Q_{10},	&& \left[Q_{5}, Q_{7}\right] = -Q_{2} + \sqrt{ 5/2}   Q_{9},\\ 
\left[Q_{5}, Q_{8}\right] = Q_{4},		&&	\left[Q_{5}, Q_{9}\right] = 2\sqrt{ 3/2}   Q_{1} - \sqrt{ 5/2}  Q_{7},\\ 
\left[Q_{5}, Q_{10}\right] = \sqrt{ 5/2}   Q_{6},	&&	\left[Q_{6}, Q_{7}\right] = - Q_{1} + 2 Q_{8},\\ 
\left[Q_{6}, Q_{8}\right] = - 2 Q_{7},	&&	\left[Q_{6}, Q_{9}\right] = \sqrt{ 3/2}   Q_{2} + \sqrt{ 5/2}  Q_{4},\\ 
\left[Q_{6},Q_{10}\right]=\sqrt{ 3/2}  Q_{3} - \sqrt{ 5/2}  Q_{5},  && \left[Q_{7}, Q_{8}\right] = 2 Q_{6}\\ 
\left[Q_{7},Q_{9}\right]=\sqrt{ 3/2} Q_{3} +\sqrt{ 5/2}  Q_{5}, &&\left[Q_{7}, Q_{10}\right] = -\sqrt{ 3/2}  Q_{2} + \sqrt{ 5/2}  Q_{4},\\ 
\left[Q_{8}, Q_{9}\right] =Q_{10},		&&	\left[Q_{8}, Q_{10}\right] = - Q_{9},\\ 
\left[Q_{9}, Q_{10}\right] = Q_{8}.\\ 
\end{array}
\end{equation}	

In order to be able to determine
the explicit form of the Jacobi operator along an arbitrary 
geodesic
$\gamma$ with initial vector $v$ at the origin $o$, it is useful to
determine the values of $ R\sp{i)}_0$, $i= 0, 1, 2, 3, 4$.
 In the following we always 
suppose that $v \in {\bf m}$ is given by 
 $$v=
\sum\nolimits_1^7 {x_i Q_i },\quad  \sum\nolimits_1^7 {(x_i )^2 }  =
1.$$  We denote $\{E_i, i = 1, \dots, 7\}$ 
the orthonormal frame field along $\gamma$ obtained by parallel 
translation of the
basis $\{Q_i \}$ along $\gamma$.

For the manifold $V_1$ the operators 
$R\sp{i)}_0,i = 0,
1, 2, 3, 4,$ written in matrix form are given by $$R\sp{i)}_0 = \left( 
{\begin{array}{*{20}c}  
{R^{i)} _{11} } & \cdots &  {R^{i)} _{17} }  \\   \cdot& 
 
& \cdot  \\     \cdot& 
 
& \cdot  \\ 
{R^{i)} _{71} } & \cdots & {R^{i)} _{77} }  \\ \end{array} }
\right)(0)$$
where ${R^{i)}}_{jk}(0)= <{R^{i)}}(E_k), E_j>(0)$.

In \cite{[T]}  Tsukada defines the curves of
constant osculator rank in the Euclidean space  and this 
concept is
applied to naturally reductive homogeneous spaces; see also 
\cite [Vol.IV,
Ch.7, Add. 4]{[Sp]}. For a unit vector $v \in {\bf m}$ 
determining the
geodesic  $\gamma$, $R_{t} = e^{t\Lambda(v)} R_0$ is a curve in 
$R({\bf m})$.
Since $e^{t\Lambda(v)}$ is a $1$-parameter subgroup of 
the group of
linear isometries of $R({\bf m})$, the curve $R_t$ has constant 
osculating rank
$r$ \cite{[T]}. Therefore, for the Jacobi operator we have
$$R_t = R_0 + a_1(t)
R\sp{1)}_0 + \cdots + a_r(t)R^{r)}_0.$$

With the help of 
 Propositions 2.3 and 2.4 we obtain:
		
\begin{lema} At $\gamma(0)$ we have:		
\begin{enumerate}
\item[i)]
$R\sp{3)}_0= -{\vert \vert \gamma'\vert \vert }^2 R\sp{1)}_0 = -
R\sp{1)}_0$;		
\item[ii)] $ R\sp{4)}_0= -{\vert \vert \gamma'\vert \vert}^2 
R\sp{2)}_0
= - R\sp{2)}_0$.		
\end{enumerate}
\end{lema}

{\bf{ Proof}}. Due to Tsukada's result about the constant
osculator rank of the curvature operator on naturally reductive spaces, we know that there exists 
$r {\in} {\Bbb N}$ such that $R^{1)}, \dots , R^{r+1)}$ are linearly dependent.
Now we are going to prove that $r = 2$ in $V_1$. 

For that we study the relationship between $R^{1)}$ and $R^{3)}$ (later we shall find another one between  $R^{2)}$ and $R^{4)}$ and so on). In particular, we have to compare $R^{1)}_{(i,j)}$ and $R^{3)}_{(i,j)}$ for $i,j=1,\dots,7$. Let us show how to proceed, for instance,  to make the comparison between $R^{1)}_{(1,1)}$ and $R^{3)}_{(1,1)}$. The computation on the other $48$ elements of $R^{1)}$ and $R^{3)} $ will be analogous.

From  Proposition 2.4 we have
\begin {eqnarray}\label{(8.1)} 
R\sp{1)}_0(X)&=& (1/2) ( [[[X,v]_{\bf{h}},v]_{\bf{m}},v]_{\bf{m}}  - [[[X,v]_{\bf{m}},v]_{\bf{h}},v]_{\bf{m}})\\ \nonumber
&=&   (1/2)     \sum\limits_{1 \leq i,j,k  \leq 7}  x_i x_j x_k ( [[[X, Q_i]_{\bf{h}}, Q_j]_{\bf{m}}, Q_k]_{\bf{m}}  - [[[X,Q_i]_{\bf{m}},Q_j]_{\bf{h}}, Q_k]_{\bf{m}}). 
  \end {eqnarray}
  Therefore if we denote $$T_1 [1,i,j,k,1] = <( \nabla_{Q_k}R)(Q_1, Q_i)Q_j, Q_1>,$$
  putting  $X= Q_1$ and using (\ref {(8.1)}) it follows
  
\begin {eqnarray}  \label{nose}
R^{1)}_{(1,1)} &=&
 <R\sp{1)}_0(Q_1), Q_1> =  (1/2) < [[[Q_1,v]_{\bf{h}},v]_{\bf{m}},v]_{\bf{m}}  
   - [[[Q_1,v]_{\bf{m}},v]_{\bf{h}},v]_{\bf{m}}, Q_1> \\ \nonumber
  &=& 
 (1/2)     \sum\limits_{1 \leq i,j,k  \leq 7}  x_i x_j x_k < [[[Q_1, Q_i]_{\bf{h}}, Q_j]_{\bf{m}}, Q_k]_{\bf{m}}  - [[[Q_1,Q_i]_{\bf{m}},Q_j]_{\bf{h}}, Q_k]_{\bf{m}}, Q_1> \\ \nonumber
 &=&  (1/2)     \sum\limits_{1 \leq i,j,k  \leq 7}  x_i x_j x_k T_1 [1,i,j,k,1]. 
   \end {eqnarray}
  
Now, using the values of the brackets of the vectors $Q_i$ in (\ref{CORCHETE}) we obtain that the  non-vanishing components of $T_1$ are 
  
\begin{equation} {\label{T1}}
\begin{array}{lclcl}
 T_1[1,2,6,4,1] = -3/2, && T_1[1,2,7,5,1] = 3/2, && T_1[1,3,6,5,1] = -3/2 \\	 
 T_1[1,3,7,4,1] = -3/2,&& T_1[1,4,2,6,1]= 3/2,&&T_1[1,4,3,7,1] = -3/2\\
T_1[1,4,4,6,1]= -\sqrt{15}/2,&& T_1[1,4,5,7,1]= -\sqrt{15}/2, && T_1[1,4,6,2,1] = -3/2\\
 T_1[1,4,6,4,1] =  -\sqrt{15},&&T_1[1,4,7,3,1] = -3/2,&& T_1[1,4,7,5,1] =  -\sqrt{15}\\
 T_1[1,5,2,7,1]= -3/2, &&T_1[1,5,3,6,1] = -3/2, &&T_1[1,5,4,7,1] = -\sqrt{15}/2\\
  T_1[1,5,5,6,1] = \sqrt{15}/2,&& T_1[1,5,6,3,1] = -3/2,&& T_1[1,5,6,5,1] = \sqrt{15}\\
   T_1[1,5,7,2,1] = 3/2&&   T_1[1,5,7,4,1] = -\sqrt{15}, &&T_1[1,6,2,4,1] = 3/2\\
  T_1[1,6,3,5,1] = 3/2 &&  T_1[1,6,4,4,1] = -\sqrt{15}/2, && T_1[1,6,5,5,1] = -\sqrt{15}/2\\
 T_1[1,7,2,5,1] = -3/2 &&  T_1[1,7,3,4,1] = 3/2,&& T_1[1,7,4,5,1]= -\sqrt{15}/2\\
   T_1[1,7,5,4,1] = -\sqrt{15}/2 \\
 \end{array}
 \end{equation}
 
 Finally using (\ref{nose}) and (\ref{T1}) it is a straightforward computation to obtain that
 
 \begin{equation}\label{R11}
 R^{1)}_{(1,1)} =   \sum\limits_{1\leq i,j,k \leq7} x_i x_j x_k T_1 [1,i,j,k,1] = -2\sqrt{15}({x_4}^2 x_6 - {x_5}^2 x_6 + 2{x_4}x_5 x_7).
  \end{equation} 
  
 For $R^{3)}$, in an analogous way, we have
\begin{eqnarray} \label{(9.1)}
R\sp{3)}_0(X)&=&   (1/8) ([[[[[X,v]_{\bf{h}},v]_{\bf{m}},v]_{\bf{m}}, v]_{\bf{m}},v]_{\bf{m}} - 3[[[[[X,v]_{\bf{m}},v]_{\bf{h}},v]_{\bf{m}}, v]_{\bf{m}},v]_{\bf{m}}  +\\ \nonumber
&&+ 3 [[[[[X,v]_{{\bf{m}}},v]_{\bf{m}},v]_{\bf{h}}, v]_{\bf{m}},v]_{\bf{m}} - [[[[[X,v]_{{\bf{m}}},v]_{\bf{m}},v]_{\bf{m}},v]_{\bf{h}},v]_{\bf{m}}). 
\end{eqnarray}

Let  $R^{3)}_{(1,1)} $ be the element  $<R\sp{3)}_0(Q_1), Q_1>$ of the matrix of $R^{3)}$. We denote 
  $$T_3 [1,i,j,k,l,m,1] = \,\, <(  \nabla\sb{ Q_m}\nabla\sb{Q_l} \nabla\sb{Q_k}R)(Q_1,Q_i)Q_j, Q_1>, \quad  i,j ,k,l,m = 1,\dots,7.$$ 
  
 Firstly, we compute the values  
   $T_3[1,i,j,k,l,m,1] $  and we compare them with the values we have obtained before for $T_1$. So, for example,  if we study the values of $T_3$ when  $i,j,k,l,m$ are respectively $1,1,4,4,6$, or one permutation $\sigma $ of these values, we obtain 
    $$ \sum\limits_{(i,j,k,l,m) \in S(1,1,4,4,6)}  x_i x_j x_k x_l x_m T_3 [1,i,j,k,l,m,1]=   (2\sqrt{15}) {x_1}^2 {x_4}^2 x_6.$$
   
 Analogously, if we consider that $i=1$, $j=1$, $k=5$, $l=5$, $m=6$, or one permutation of these values,  we obtain 
    $$ \sum\limits_{(i,j,k,l,m) \in S(1,1,5,5,6)}  x_i x_j x_k x_l x_m T_3 [1,i,j,k,l,m,1]=   (-2\sqrt{15}) {x_1}^2 {x_5}^2 x_6$$
    and for $i=1$, $j=1$, $k=4$, $l=5$, $m=7$,  or one permutation of these values, it has
       $$ \sum\limits_{(i,j,k,l,m) \in S(1,1,4,5,7)}  x_i x_j x_k x_l x_m T_3 [1,i,j,k,l,m,1]=   (4\sqrt{15}) {x_1}^2 x_4 {x_5} x_7.$$

  On the other hand, we have that for other sets of indices $I$ different from $A$, $B$ or $C$  where
  $A=\{h,h,4,4,6\}$,  $B=\{h,h,5,5,6\}$,  $C=\{h,h,4,5,7\}$, $h= 1,\dots,7$,  the following sum vanishes
    $$ \sum\limits_{(i,j,k,l,m) \in S(I)}  x_i x_j x_k x_l x_m T_3 [1,i,j,k,l,m,1]=  0.$$
 In consequence
\begin{equation}\label{R3}
 \sum\limits_{1 \leq i,j, k,l,m  \leq 7}  x_i x_j x_k x_l x_m T_3 [1,i,j,k,l,m,1]=   \sum\limits_{1 \leq h \leq 7} ({x_h}^2) \, 2\sqrt{15}({x_4}^2 x_6 - {x_5}^2 x_6 + 2{x_4}x_5 x_7)  .
  \end{equation}

Then  we  can conclude from (\ref{R11}) and (\ref{R3}) that
 $$R^{3)}_{(1,1)}  = - ({x_1}^2 + \cdots + {x_7}^2) R^{1)}_{(1,1)}  = -{\vert \vert \gamma'\vert \vert }^2 R^{1)}_{(1,1)}.$$
 
 The proof of ii) is analogous to i).
 
Remark .  Using {\it{Mathematica}}  it is possible to prove that the non-null components we have calculated in the proof are correct.
 
\begin{prop}			At $\gamma(0)$ we have:		
\begin{enumerate}
\item[i)]	
$R\sp{2n)}_0 =
(-1)\sp{n-1} R\sp{2)}_0$;
\item[ii)] $R\sp{2n+1)}_0 = (-1)\sp{n}
R\sp{1)}_0$.
\end{enumerate}

\end{prop}	
\pf		We are going to prove i) by induction, ii) may be obtained in a 
similar way. First, Lemma
3.1  (ii) gives the result for $n = 2$. Next, suppose that for $n = k$ the result 
is true,
that is,	
\begin {equation}\label{(3.1)}
 R\sp{2k)}_0 = (-1)^{k-1} R^{2)}_{0}.	
\end {equation}
	Using
Proposition 2.4 we have
$$(-1)\sp{2k+1} 2\sp{2k+2 }R\sp{2k+2)}_0(X) =
\sum\limits_{i = 0}^{2k + 2} 		{( - 1)^i } \left( 
{\begin{array}{*{20}c}  
{2k + 2}  \\   i  \\ \end{array} } \right)[[[X,v]_{\bf 
m},\dots, v]_ {\bf h}^{i+1)}, \dots
, v]_{{\bf m}}^{2k+4)}.$$			
	There are $2k+3$ terms and each one has $2k+4$
brackets. If we take into account Lemma 1.1 in the previous expression, we
obtain
\begin{eqnarray*}	&&	(-1)^{2k+1} (2)^{2k+2} R\sp{2k+2)}_0(X)
\\	&=&\sum\limits_{i = 0}^{2k} {( - 1)^i } \left( 
{\begin{array}{*{20}c}  
{2k}  \\   i  \\ \end{array} } \right) [[[[[X,v]_{\bf 
m}, \dots,
v]_{\bf h} ^{ i+1)}, \dots 		,v]_{\bf m}^{ 
2k+2)},v]_{\bf 
m},v]_{\bf
m}\\		&& -2 (\sum\limits_{i 		= 0}^{2k} {( - 1)^i }
\left({\begin{array}{*{20}c}   {2k}  \\   i  \\ \end{array} } 
\right)
 [[[[X,v]_{\bf m},\dots , v]_{\bf h}^{ 
i+2)},\dots
, v]_{{\bf m}} ^{ 2k+3)},v]_{\bf m})\\		&& +\sum\limits_{i = 
0}^{2k} {( - 
1)^i
} \left( {\begin{array}{*{20}c}   {2k}  \\   i  \\\end{array} } 
\right)
[[[X,v]_{\bf m},\dots , v]_{\bf h}^{ i+3)},\dots
, v]_{\bf m}^{ 2k+4)}.\end{eqnarray*}	

If we call 		$X' =[X,v]_{\bf m}$ and $X''
= [[X,v]_{\bf m},v]_{\bf m}$ and using also Proposition 2.4, it
follows
\begin{eqnarray*}		&&(-1)^{2k+1} 2^{2k+2} R\sp{2k+2)}_0(X)\\
		&=&
(-1)^{2k+1} 2^{2k+2} ([[R\sp{2k)}_0(X),v]_{\bf m},v]_{\bf m} -
2[R\sp{2k)}_0(X'),v]_{\bf m} +R\sp{2k)}_0(X'')).		
\end{eqnarray*}
Taking into account Lemma 3.1, formula (\ref{(3.1)}) and the values of $X'$ 
and $X''$ in
previous expression, Proposition 3.2 follows.

The next result follows immediately
from Proposition 3.2.				
\begin{prop}			 The normal
naturally reductive homogeneous space $V_1 = Sp(2)/SU(2)$ is of 
constant
osculator rank $2$.			
\end{prop}		
\begin{cor}			
 Along the geodesic $\gamma$ the Jacobi
operator can be written as 
$$R_t= R_0 + R\sp{2)}_0 + R\sp{1)}_0\sin t - 
R\sp{2)}_0
\cos t.$$			
\end{cor}	
The proof  follows from the Taylor development of $R_t$ at $t=0$ and by using
Proposition 3.2.

\begin{cor}				Along
the geodesic $\gamma$ the derivatives of the Jacobi operator
satisfy:		
\begin{enumerate}
\item[i)]
$R\sp{2n)}_t  = (-1)\sp{n-1} R\sp{2)}\sb{t}$;		
\item[ii)]
$R\sp{2n+1)}\sb{t}  = (-1)\sp{n} R\sp{1)}\sb{t}$;		
\item[iii)]
 $ R_t \cdot R\sp{1)}_t= R\sp{1)}_t \cdot R_t ,\\ R_t \cdot R\sp{2)}_t = 
R\sp{2)}_t \cdot R_t,
\\R\sp{1)}_t \cdot R\sp{2)}_t= R\sp{2)}_t \cdot R\sp{1)}\sb{t}$.
\end{enumerate}
\end{cor}
The result is a consequence of Corollary 3.4 and  the fact that iii) is true for $t=0$.
\begin {rem} {\rm In \cite{[BV]} the authors analyze a class of Riemannian 
homogeneous 
spaces
 which have the property that the eigenspaces of $R_{\gamma}$ are 
parallel along
$\gamma$. Evidently, this property is not verified in our case. In 
fact, although,
for each $t$, the operators $R_t$ and $R\sp{i)}_t $ commute and 
therefore are simultaneously diagonalizable, according to Corollary 3.4 and (ii) 
of Proposition 2.2 the eigenvectors of these operators are not independent of $t$.}
\end {rem}

\section{The solution of the
Jacobi equation on the manifold $V_1$. Application to the determination of volumes of
geodesic balls.} For a naturally  reductive
homogeneous Riemannian manifold it is possible to write the Jacobi 
equation
as a differential equation with constant coefficients. In order to 
do that,
the canonical connection is frequently used. Since this connection 
and the
Levi-Civita connection have the same geodesics, in an equivalent 
form, it is
possible to write the same equation based on the Levi-Civita 
connection 
\cite{[Ch2]}.
In this case, the coefficients are functions of the arc-length along 
the
geodesic. In order to work with this equation on the manifold $V_1$ it 
will be
useful to use the simple expression of the Jacobi operator $R_t$.	 
	We shall
now introduce some notation and provide some basic formulae which 
will be
needed in this section. For more information see 
\cite{[BPV],[Ch2],[Z]}.	
Let $A$ be the Jacobi tensor field along the geodesic $\gamma$ (that 
is,
the solution of the endomorphism valued Jacobi equation $Y'' + R_t 
Y= 0$ along
$\gamma$) with initial values		
\begin {equation}\label{(4.1)}
A_0 = 0,\quad  A\sp{1)}_0 =I,
\end {equation}					
		 where we consider the covariant differentiation with respect
to $\gamma'$ and $I$ is the identity transformation of 
$T_{\gamma(0)}M$. Then,
the Jacobi's equation is $A^{2)}_t= - A\sb{t} R_t $. 

In order to be able to obtain the expression 
of the
Jacobi fields with initial conditions (\ref{(4.1)}) at $\gamma(0)$, it is 
enough to know
the development in Taylor's series of $A\sb{t}$ and to apply the 
initial
conditions. Thus, using Lemma 3.1, in the power series of 
$A\sb{t}$ only 
appear $R_0$,
$R\sp{1)}_0$ and $R\sp{2)}_0$.	If $\{E_{i}, i=1,\dots,7 \}$ 	is the 
orthonormal frame 
field
along $\gamma$ obtained by parallel translation of the basis $\{Q_i\}$
along $\gamma$, one has $Y\sb{t} = A\sb{t} E\sb{t}$	or
$Y\sb{i,t} = A^{j}\sb{i,t}E_{j,t}$ , $1 \le i, j \le 7$, and this is 
the 
expression
of the Jacobi vector fields along the geodesic $\gamma$ with 
the indicated
initial conditions.  
\begin{prop}
For the manifold $V_1$, one has
			$$A\sb{t}=  	\sum\limits_{k =0}^\infty  \frac{{{1}}}{{{{k!}} }} {\beta  _k  t^k  } $$				
where
$\beta_k =
\alpha_{k-1}+\beta_{k-1}'$, \quad $\alpha_k = \alpha_{k-1}'-R\beta_{k-1}$, $k \geq  2$.
Moreover, $\alpha_0  = \beta_0 = 0$, 		$\alpha_1  = 0$,
	$\beta_1 = 
I$			and the
coefficients 			$\beta_k$ are only functions of $R_0$, 
$R^{1)}_0$ and
$R^{2)}_0$.						
\end{prop}
\pf If we successively derive $A_t^{2)} = -R_tA_t$, we have
$$A_t^{i) }= \left(\alpha'_{i-1}(t)-R_t\beta_{i-1}(t)\right)A_t + \left(\alpha_{i-
1}(t)+\beta'_{i-
1}(t)
\right)A_t^{1)},$$  we can write this expression as
	$$A_t^{i)} = \alpha_i(t) A_t + \beta_i(t) A_t^{1)},$$ 
where 
$$\alpha_i(t) =
\alpha_{i-1}'(t)-R_t\beta_{i-1}(t),$$ and $$ \beta_i(t) = \alpha_{i-
1}(t)+\beta_{i-1}'(t),
 \quad i \geq
2;$$
if $t = 0$ one has $A_0^{0)}= \beta_0(0)= 0$, 	$A^{1)}_0 = \beta_1(0) 
=I$,	
$A^{2)}_0= \beta _2 (0)= 0$,	$A^{3)}_0= \beta_3 (0)= -R_0^{0)}=-R_0$, 
and, in 
general,
$$A^{i)}_0 = \alpha_{i-1}(0)+\beta_{i-1}'(0) = \beta_i(0).$$

 If there 
is
no confusion we will identify $ \alpha_{i} = \alpha_{i}(0)$ and 
$ \beta_{i} = \beta_{i}(0)$.
Now the result follows using the development in Taylor's series of $A_t$.

	Let $m$ be a point of the manifold $M$ and $V$
and $U$ open neighbourhoods of $0$ in $T_mM$ and of $m$ in $M$ 
respectively such that
 $\exp_m$ is a diffeomorphism of $V$ onto $U$. For all $v \in V$, 
$\theta(v)$ \cite [p. 54]{[BGM]} is a
well-defined function, it is defined as the absolute value
of a determinant function: 
$$\theta(v)=  \vert \det T_v \exp_m \vert.$$

\begin{defi}Let $U_{\epsilon}(m)$ be a normal neighbourhood of 
radius  
${\epsilon}>0$ of the point $m$ in $M$. For each $t$ such that 
$0<t<\epsilon$ 
and for each $v$ in $T_mM$ the function $t \mapsto \theta(tv)$ is 
the volume
density function at $m$ in the direction $v$. 
\end{defi}
\begin{lema}{\em \cite [p. 90]{[BGM]}} 
Let $u \in T_o M$ and $t>0$, then for
all $v \in T_o M$, $T\sb{tu}\exp_o(v)$ is the value in $t$ of the 
Jacobi field
$Y$ along the geodesic $\gamma$ ($\gamma(0)=o, \gamma'(0)=u$) with 
initial
conditions $Y(o) = 0$, $Y'(o) = v/t$.\end{lema}
\begin{prop}					In the
manifold $V_1$, the volume density function at $o$ is given by
\begin{equation}
\theta(tu) = 
\frac{{{1}}}{{{{t}}^{{7}} }}\left| {\det A} \right|.
\end{equation}    
\end{prop}	
The proof follows in a natural way from the standard methods of 
\cite{[BGM],[Gr], 
[GV]}.
\begin{cor}				The coefficient of $t\sp{n}$ in the
development of $\det A$ is given by
\[ a_n={{ }} \sum\limits_{
\matrix{  r_1 +
\cdots + r_7 = n \cr 
0\leq r_1,\cdots,r_7 \leq n\\ }} \frac{{{1}}}{{{{r_{1}!}} }} \cdots \frac{{{1}}}{{{{r_{7}!}}
}}
{\sum\limits_\sigma  {\sig(\sigma )} \beta _{r_1, \sigma  (1)}^{1} \cdots\beta
_{r_7, \sigma (7)}^{7} }.
\]\end{cor}								
\pf				The seven columns $C_j, j=1,\dots,7$ of $A$ can
be written as
\[C_j  = \beta _0 ^j  +\cdots + \frac{{{1}}}{{{{n!}} }} \beta _n ^j t^n 
+ \cdots =  \sum\limits_{k = 0}^\infty  \frac{{{1}}}{{{{k!}} }} {\beta _k^j } t^k \]
where the upper index
$j$ shows the $j^{th}$-column of the matrix $\beta_ k$ (Proposition 4.1). 
Taking into
account that the determinant is a multilinear function, the 
coefficient $a_n$
of $t\sp{n}$ in the development of $\det A$ is  
\[ {{a}}_{{n}} {{  =
}}  \sum\limits_{\matrix{  r_1 +
\cdots + r_7 = n \cr 
0\leq r_1,\dots,r_7 \leq n\\ }} \frac{{{1}}}{{{{r_{1}!}} }} \cdots \frac{{{1}}}{{{{r_{7}!}} }} {\det (\beta _{r_1}^1 
,\dots,\beta _{r_7}^7  ).}
\] If we represent the matrix $ \beta_{r_{k}} $ by $\beta_{r_{k}} = (\beta_{r_{k},i}^{j} ), i, j
= 1,\dots , 7$,  using
the algebraic definition of the determinant it follows that				
$$ a_n={{ }} \sum\limits_{\matrix{  r_1 +
\cdots + r_7 = n \cr
0\leq r_1,\dots,r_7 \leq n\\ }} \frac{{{1}}}{{{{r_{1}!}} }} \cdots \frac{{{1}}}{{{{r_{7}!}}
}}
{\sum\limits_\sigma  {\sig(\sigma )} \beta _{r_1, \sigma  (1)}^{1} \cdots\beta
_{r_7, \sigma (7)}^{7} }$$
 where $\sigma$ are the
permutations of seven elements and $\sig(\sigma)$ represents the 
signature of
the corresponding permutation. 							
\begin{lema}			{\em \cite{[Gr]}} For the manifold
$V_1$, 		
\begin{enumerate}
\item[i)] The area of the geodesic sphere with center $o \in V_1$ and 
radius
$t$ is given by $$S_o(t)= t^6 \int\nolimits_{\Omega ^6 (1)} {\theta 
(tu)du} $$ 
 where 
$\Omega^6(1)$ denotes the 6-dimensional Euclidean unit sphere.

\item[ii)] The volume of the geodesic ball with center $o\in V_1$ 
and
radius $r$ is given by
$$V_o(r)=\int\nolimits_0^r {S_o(t)dt}.$$
\end{enumerate}
\end{lema}						
Now, using the standard notation for moments \cite [p. 255--258]{[Gr]}, we have:
\begin{prop}			
\begin{enumerate}
\item[i)] The area of the geodesic ball with
center $o$ and radius $t$ is given by
$$S_o(t) = (16\pi^{3}/105 )  {{{
}}}\sum\limits_{n = 3}^\infty  < a_{2n+1} > t^{2n};$$
\item[ii)]
The volume of the geodesic ball with center $o$ and radius $t$ is 
given
by			$$ V_o(t)=(16\pi^{3}/105)  
\sum\limits_{n = 3}^\infty \frac{{{1}}}{{{{2n+1}} }} < a_{2n+1} > t^{2n+1}.$$			
\end{enumerate}
\end{prop}				
\pf			i) If we integrate i) of Lemma 4.6 over the sphere we 
have that the odd 
powers vanish and then the result follows immediately. For ii) we use that
$V_o(r)=\int\nolimits_0^r  {S_o}(t)dt$.\par

{\bf Acknowledgements}: The authors gladly acknowledge helpful
conversations with  J. \'Alvarez Paiva, T. Arias--Marco, J. C. Gonz\'alez -
D\'avila, 
O. Kowalski, E. Mac\'{\i}as
and L. Vanhecke.			

\newpage
Author's adresses:

A. M. Naveira\\
Departamento de Geometr\'{\i}a y Topolog\'{\i}a. Facultad de
Matem\'aticas.\\
Avda. Andr\'es Estell\'es, N¼ 1\\
46100 - Burjassot\\
Valencia,
SPAIN\par
Phone  +34-963544363\\
Fax: 
+34-963544571\\
e-mail:	{\tt naveira@uv.es}

A. D. Tarr\'{\i}o Tobar\\
E. U.
Arquitectura T\'ecnica\\
Campus A Zapateira. Universidad de  A Coru\~na\\
15192 - A
Coru\~na, SPAIN\par
Phone  +34-981167000 Ext. 2721, 2713\\
Fax: +34-981167060\\
e-mail:
{\tt madorana@udc.es}
\end{document}